\documentclass{article}
\usepackage{amsfonts,amssymb,latexsym}
\newcounter{num}[section]
\newcommand{\Num}{\refstepcounter{num}%
\textbf{\arabic{section}.\arabic{num}}}

\newcommand{\Theorem}{\textbf{Theorem~}}
\newcommand{\Proof}{{\sc{Proof}}}
\newcommand{\Def}{\textbf{Definition~}}

\newcommand{\Lemma}{ \textbf{Lemma~}}
\newcommand{\Ex}{  \textbf{Example}}

\newcommand{\Prop}{\textbf{Proposition~}}
\newcommand{\Cor}{ \textbf{ Corollary~}}

\newcommand{\Ac}{{\cal A}}
\newcommand{\Bc}{{\cal B}}

\newcommand{\Nc}{{\mathcal N}}
\newcommand{\Lc}{{\mathcal L}}
\newcommand{\Dc}{{\cal D}}

\newcommand{\Ic}{\mathcal{I}}

\newcommand{\Ec}{\mathcal{E}}
\newcommand{\Fc}{\mathcal{F}}
\newcommand{\al}{{\alpha}}
\newcommand{\la}{{\lambda}}
\newcommand{\La}{{\Lambda}}
\newcommand{\Uc}{{\mathcal{U}}}
\newcommand{\Sc}{{\mathcal{S}}}
\newcommand{\Rc}{{\mathcal{R}}}
\newcommand{\Xc}{{\mathcal{X}}}

\newcommand{\Fq}{{\Bbb F}_q}
\newcommand{\Cb}{{\Bbb C}}

\newcommand{\Lb}{{\mathbb{L}}}
\newcommand{\Ub}{{\mathbb{U}}}
\newcommand{\Ib}{{\mathbb{I}}}
\newcommand{\Jb}{{\mathbb{J}}}
\newcommand{\Pb}{{\mathbb{P}}}
\newcommand{\Nb}{{\mathbb{N}}}
\newcommand{\Bb}{{\mathbb{B}}}
\newcommand{\Gb}{{\mathbb{G}}}
\newcommand{\Hb}{{\mathbb{H}}}

\newcommand{\Add}{{\mathrm{Ad}}}

\newcommand{\GL}{{\mathrm{GL}}}

\newcommand{\Ind}{{\mathrm{Ind}}}
\newcommand{\Or}{{\mathrm{O}}}
\newcommand{\Sp}{{\mathrm{Sp}}}
\newcommand{\UT}{{\mathrm{UT}}}

\newcommand{\gx}{{\mathfrak g}}
\newcommand{\gl}{\mathfrak{gl}}
\newcommand{\nx}{{\mathfrak n}}

\newcommand{\bx}{\mathfrak{b}}
\newcommand{\ux}{{\mathfrak u}}
\newcommand{\ax}{{\mathfrak a}}
\newcommand{\Ax}{{\mathfrak{A}}}
\newcommand{\Bx}{{\mathfrak{B}}}

\newcommand{\px}{{\mathfrak p}}
\newcommand{\lx}{{\mathfrak l}}

\newcommand{\Gx}{{\mathfrak G}}

\newcommand{\eps}{{\varepsilon}}

\newcommand{\row}{{\mathrm{row}}}
\newcommand{\col}{{\mathrm{col}}}
\newcommand{\Ad}{{\mathrm{Ad}}}
\newcommand{\spann}{\mathrm{span}}

\renewcommand{\leq}{\leqslant}
\renewcommand{\geq}{\geqslant}

\begin{document}
\Large

\title{Supercharacters for parabolic contractions of finite groups of $A,B,C,D$ Lie types }
\author{A.N.\,Panov \\
	\\
\emph{Mathematical Department, Samara  University, Samara, Russia,
apanov@list.ru}}
\date{}
 \maketitle
 
 {\small
 	We construct  supercharacter theories for the  finite groups constructed by parabolic constraction from simple groups of $A,B,C,D$ Lie types. In terms of rook placements in the root systems
 	we classify  supecharacters and superclasses. \\
 	\\
 	{\bf \hspace{0.3cm} Key Words:} group representations, supercharacter theory, group contraction, irreducible character, superclasses\\
 	\\
 	{\bf \hspace{0.3cm} 2010 Mathematics Subject Classification:} 20C15, 17B22}

\section{Introduction} 
The notion of a supercharacter theory was introduced by P. Diaconis and I.M. Isaacs in  the paper  \cite{DI} in 2008.  
By definition, a supercharacter theory for a given finite group $G$ is a pair  $({\mathfrak S},{\cal K})$, where  ${\mathfrak S}=\{\chi_1,\ldots,\chi_M\}$ is a system of pairwise disjoint  characters (representations) of  $G$  and ${\cal K} =\{K_1,\ldots, K_M\}$ is a partition of $G$ such that 
the characters from  ${\mathfrak S}$ are constant on subsets from   ${\cal K}$ and $\{1\}\in {\cal K}$.
The characters from  ${\mathfrak S}$ are called  \emph{supercharacters}, and subsets from  ${\cal K}$ -- \emph{superclasses} (see \cite{DI}). Observe the number of supercharacters equals to the number of superclasses.  
An example of supercharacter theory is the  pair that consists of the system of irreducible characters and the decomposition of group into conjugacy classes. 
If a classification of irreducible representations is considered  an extremely difficult, "wild"\, problem, we aim  to construct a construc\-tive supercharacter theory.
The classical example is the theory of basic characters for the unitriangular group  $\mathrm{UT}(n,\Fq)$  developed by C.\,Andr\'{e} in  \cite{A1,A2,A3}.  This example was generalized to algebra groups in the same paper \cite{DI} (by definition, an algebra group is a group of the form  $\Gx=1+J$, where $J$ is a finite dimensional associative nilpotent algebra over a finite field).

In the present paper, we construct a supercharacter theory for parabolic contractions of the groups $\GL(n)$, ~$\Or(n)$ and $\Sp(n)$ defined over a finite field $\Fq$. For each group we classify supercharacters and superclasses.  

To construct a supercharacter theory for the parabolic contraction   $\GL^a(n)$, we use the method  of  \cite{P1} (see Theorem \ref{GGsuper}).  The case of Borel contraction  is studed in  \cite{P2}, where we apply the other method  from \cite{P3,P4,P5} that gives rise a  supercharacter theory  that is slightly different from the one constructed in the present paper. Also, we  classify supercharacters and superclasses for the unipotent  radical  $U^a$  (see Theorems  \ref{supclass} and \ref{ThSupUa}). For the  parabolic radicals of $\GL(n)$, the similar classification is verified in  \cite{Thiem}. 

A supercharacter theory for the Sylow subgroups in the orthogonal and symplec\-tic groups is constructed in  \cite{AN-1,AN-2,AFN,SA}. 
The author proposed  a super\-cha\-rac\-ter theory for parabolic subgroups in  $\Or(n)$ and $\Sp(n)$ in the paper  \cite{P6}. 
This approach is developing in the present paper for parabolic contractions of orthogonal and symplectic groups. We obtain a description of supercharacters and superclasses in the exact form (see  Theorems \ref{SupUaos} and \ref{Gsuper}).

\section{Parabolic contractions}

Let  $G$ be a split reductive groups defined over a field  $K$ and  $\gx$ be its Lie algebra. 
 Fix the Cartan subgroup $H$ contained in the Borel subgroup  $B$. The group  $B$ has the unipotent radical  $N_+$ with the Lie algebra  $\nx_+$.  Applying the Killing form  $(\cdot,\cdot)$  of  $\gx$, we identify  $\nx^*$ with the opposite subalgebra  $\nx_-$. 
 
  Let $P$ be a parabolic subgroup containing $B$. The subgroup  $P$ is a semidirect product   $P=LU_+$, where  $L$ is the Levi subgroup containing  $H$,  and  $U_+$ is the unipotent radical contained in   $ N_+$.  The Lie algebra  $\px$ of the group  $P$ is a  semidirect sum  $\px=\lx\oplus\ux_+$ of the  Lie subalgebra $\lx=\mathrm{Lie}(L)$ and  the ideal  $\ux_+=\mathrm{Lie}(U_+)$. 
 
As above we identify  $\ux_-$ with $\ux_+^*$. 
Consider the commutative Lie algebra $\ux_-^a$ that coincides with  $\ux_-$ as a linear space and the abelian group  $U_-^a= 1+\ux_-^a$  with multiplication  $(1+\la_1)(1+\la_2)=1+\la_1+\la_2$.

 \emph{A parabolic contraction} of the group  $G$ is a semidirect product   $G^a=P\ltimes U_-^a$;
 the subgroup  $P$ acts on $U^a_-$ by the formula  $p(1+\la)p^{-1} = 1+ \Add_p^*(\la)$, where $p\in P$ and $\la\in\ux_-^a=\ux^*$.  
 Respectively, its Lie algebra  $\gx^a= \px\oplus\ux_-^a$  (\emph{a parabolic contraction} of the Lie algebra $\gx$) is a semidirect sum of the Lie subalgebra $\px$ and the ideal $\ux_-^a$.
These definitions are matched with the general definition of Lie algebra contraction from \cite{IW, GOV}.

\section{Supercharacters and superclasses of the group  
$\GL^a(n)$}

Let  $\GL(n)$ be the general linear group over a finite field $\Fq$ of  $q$ elements and of characteristic $\ne 2$. Fix a parabolic subgroup  $P$ of $\GL(n)$ and consider its parabolic contraction  $\GL^a(n)$.
   The parabolic subgroup   $P$  is defined by a decomposition of  $[1,n]$ into the system of consecutive segments $\Ic=I_1\sqcup\ldots\sqcup I_\ell$. Let $n_i=|I_i|$. Then
 $L=\GL(n_1)\times\ldots\times \GL(n_\ell)$.  
     
  Emphasize  the Lie algebra $\gl^a(n)$ of the group  $\GL^a(n)$  is an associative algebra with respect to the  multiplication for which  $\px$ is a subalgebra with matrix multiplication, $\ux_-^a$ is a subalgebra with zero multiplication and 
   \begin{equation}\label{rlaction}
  t\la(x)=\la(xt)\quad\mbox{and}\quad \la t(x)=\la(tx),
  \end{equation} 
  where $t\in \px$, ~$x\in \ux_+$, and $\la\in \ux_-^a$. 
  The group $\GL^a(n)$ is realized as the subgroup of invertible elements of the associative algebra  $\gl^a(n)$.
  
  The group  $\GL^a(n)$ is a semidirect product  $\GL^a(n)=L\ltimes U^a$, where $U^a=U_+\ltimes U_-^a$.
  The subgroup $U^a$  is an algebra subgroup  $U^a=1+\ux^a$, where
   $\ux^a=\ux_++\ux_-^a$ is the nilradical of  $\gl^a(n)$. 
   
   Our goal is to construct a supercharacter theory for the group  $\GL^a(n)$ following the general method from the paper  \cite{P1} (where it is refered to as the $GG$-supercharacter theory).

    \subsection{Supercharacter theory for  $U^a$}
 
   Because $U^a$ is an algbera subgroup, it has the supercharacter theory of P. Diaconis and I.M. Isaacs from \cite{DI}. The superclasses in algebra group  $\Gx=1+J$  is the subsets of the form  $1+\Gx x\Gx$, where $x\in J$. 
   The supercharacters $\chi_\la$  of   $\Gx$ are constructed  by the following procedure.  Fix the nontrivial  character  $t\mapsto\eps^t$ of the additive group of  field $\Fq$ with values in  $\Cb^*$.  Consider the system of representatives   $\{\la\}$ of the left-right orbits of the group  $\Gx\times \Gx$ in $J^*$.  By each  $\la$, we construct the linear character  
   $$\xi_\la(1+x)=\eps^{\la(x)}$$   of the stabilizer  $\Gx_{\la,\mathrm{right}}$  of the linear form  $\la$ with respect to the right action of  $\Gx$ on $J^*$ (see formula (\ref{rlaction})). The supercharacter  $\chi_\la$ of the algebra  group $\Gx$ is induced  from $\xi_\la$ to $\Gx$.  There is a formula for supercharacter $$\chi_\la(1+x)=\frac{|\Gx\la|}{|\Gx\la \Gx|}\sum_{\mu\in \Gx\la \Gx}\eps^{\mu(x)}.$$
   So, the superclasses and supercharacters of an algebra group  $\Gx$ are uniquely determined by  $(\Gx\times \Gx)$-orbits  in $J$ and $J^*$ respectively.  
   
   In our case,  $\Gx=U^a$ and  a classification of    $(U^a\times U^a)$-orbits in $\ux^a$ and  $(\ux^a)^*$ is considered an extremely difficult problem  (see \cite{P1,P5,Thiem}).
   Therefore, we consider  the more "coarse"\, supercharacter theory  (see Theorem \ref{ThSupUa}) in which superclasses are attached to the   $(\GL^a(n)\times \GL^a(n))$-orbits in $\ux^a$,
   and supercharac\-ters a up to a constant multiple equal to  sums of functions  
   $\eps^{\mu(x)}$, where $\mu$ runs through a $(\GL^a(n)\times \GL^a(n))$-orbit in  $(\ux^a)^*$.

    Let  $\Delta$ denote the root system of the Lie algebra  $\gl(n)$. We attach to each root  $\al=\epsilon_i-\epsilon_j$,~$1\leq i,~j\leq n$,~  $i\ne j$, the pair  $(i,j)$.  To  simplify notations, we write $\al=(i,j)$.
   We refer to  $i$ as a row number of the root  $\al$ and  to $j$ as a column number of root $\al$. Notations: $i=\row(\al)$ and  $j=\col(\al)$.
    
    A positive root  $\al\in \Delta_+$ is attached to a pair with  $i<j$. Respectively, for an negative root  $\al\in \Delta_-$, we have $i>j$. 
    
    We introduce the definition of  \emph{sum of roots}, which if different from usual one:  a root $\gamma$ is a sum  $\gamma=\al+\beta$  if $\al=(i,j)$, ~ $\beta=(j,k)$ and  $\gamma=(i,k)$.
    
   To each root $\al\in\Delta_+$, we attach the matrix unit  $E_\al\in\nx_+$ and to each root
     $\al\in \Delta_-$ -- the matrix unit  $F_\al\in\nx_-^a$.

    Consider the subset   $\Delta_+(\ux^a)$ which consists of  $\al\in \Delta_+$ such that  $E_\al\in \ux_+$. Respectively,  $\Delta_-(\ux^a)=\{\al\in\Delta_-:~ F_\al\in\ux_-^a\}$ and $\Delta(\ux^a)=\Delta_+(\ux^a)\cup\Delta_-(\ux^a)$.

    We realize $\ux^a$ as the subspace of  the linear space of  $(n\times n)$-matrices spanned by the union of two systems of matrix units $$\{E_\al:~ \al\in \Delta_+(\ux^a)\}\cup\{F_\al:~ \al\in \Delta_-(\ux^a)\}.$$
 The structure relations have the form  
  \begin{itemize}
  	\item  $E_\al E_\beta = E_{\al+\beta}$ if  $\al+\beta$ is defined, otherwise  $E_\al E_\beta=0$.
  	\item  $F_\al F_\beta =0$.
  	\item  $E_\al F_\beta = F_{\al+\beta}$ if $\al+\beta$ is defined and belongs to  $\Delta_-(\ux^a)$, otherwise   $E_\al F_\beta=0$.
  	\item $F_\al E_\beta = F_{\al+\beta}$  if  $\al+\beta$ is defined and belongs to  $\Delta_-(\ux^a)$,  otherwise   $F_\al E_\beta=0$.
  \end{itemize}
    
Observe that    $\nx_+=\ux_+\oplus \nx_+(L)$, where $\nx_+(L)=\nx_+\cap \lx$. Respectively,  $\nx_-^a= \ux_-^a\oplus \nx^a_-(L)$, where $\nx^a_-(L)=\{\la\in\nx_-^a:~ \la(\ux_+)=0\}$.\\
 \Lemma\Num\label{muux}. 1) $\nx^a_-(L)$ is a two sided ideal in the associative algebra $\gl^a(n)$.\\
 2) $\mu_Lx=x\nu_L=0$  for all  $x\in \ux_+$ and  $\mu_L,~\nu_L\in \nx^a_-(L)$. \\
 \Proof. The both statements follows from the  fact that  $\ux_+$ is a two sided ideal in  $\nx_+$. $\Box$
 
 For an arbitrary element   $u=1+X\in U^a$, where $X\in \ux^a$, we consider the class $$K(u)=\{1+AXB: ~ A,B\in \GL^a(n)\}.$$   
 The following Theorem  \ref{supclass} gives a classification of classes  $K(u)$. We need new definitions.  
  
  \emph{A rook spacement } in  $\Delta$ is a subset  $D\subset \Delta$
  that has at most one root in any row and colunm.
 The rook spacement  $D$ is a union  $D=D_+\cup D_-$, where $D_\pm=D\cap \Delta_\pm$. 
  In his papers, C.\, Andr\'{e} used the other term  \emph{basic subset} for the root placements in  $\Delta_+$.

  To each rook placement  $D\subset \Delta(\ux^a)$ and  a map  $\phi: D\to \Fq^*$, we attach the element   $X_{D,\phi}\in \ux^a$ of the form
  
  \begin{equation}\label{xdphi}
 X_{D,\phi} = \sum_{\gamma_+\in D_+}\phi(\gamma_+)E_{\gamma_+}  + 
 \sum_{\gamma_-\in D_-}\phi(\gamma_-)F_{\gamma_-}.
  \end{equation}
  
 If $\phi$ is identically equals to one on  $D$, we take $X_{D,\phi}=X_D$. 
    Denote $u_D=1+X_D$. 
    
     For any rook placement   $D\subset \Delta$ and pairs of segments  $I_k, I_m$, where $k\ne m$, we define 
   $$r_{ij} = r_{ij}(D) = \mathrm{rank}\, X_D(I_k\times I_m) = |D\cap(I_k\times I_m)|,$$
   $$n'_{km} = \sum_{s<k}n_s +\sum_{t<m}r_{kt},$$
   $$n''_{km} = \sum_{s<m}n_s +\sum_{t<k}r_{tm}.$$
  We say that a rook placement  $D_c$  has  a  \emph{canonical type} if  for each pair  $1\leq k, m\leq \ell$,~$k\ne m$, the subset  
   $D_c\cap(I_k\times I_m)$ consists of roots  
   $$(n_{km}'+1, n_{km}''+1),\ldots, (n_{km}'+r_{km}, n_{km}''+r_{km}).$$
   Let  $W_L$ be a Weyl group in $L$. For each  $w\in W_L$,  the action of $w$ on the row numbers (respectively, the action of  $w^{-1}$ on the column numbers) of roots from   $\Delta$ defines  the left (respectively, the right) action of  $w$ on $\Delta$.\\
   \Lemma\Num\label{canon}. For any rook placement  $D$ there exist  elements $w_1,w_2\in W_L$ such that  $w_1Dw_2$ is a rook placement  $D_c$ of canonical type.  The rook placement  $D_c$ is uniquely determined by  $D$.\\
   \Proof.  For any rook placement $D$, there exists  $w_1\in W_L$  such that the  system  $\row(w_1D)$ has the required ordering.  Multiplying by a suitable  $w_2\in W_L$ we obtain the  rook placement  $D_c=w_1Dw_2$ of canonical type. Uniqueness follows from the fact that the numbers $\{r_{km}\}$ don't change if we replace $D$ by $w_1Dw_2$.~$\Box$
    
  In the case  $P=B$, the  next theorem is proved in  \cite{P2}.\\  
 \Theorem\Num\label{supclass}. 1) Each class contains the element  $ u_{D}$ for some rook placement  $D$ in the rook system  $\Delta(\ux^a)$. \\
 2) Two elements  $ u_{D}$ and $ u_{D'}$ belong to a common class if and only if $D$ and $D'$ are conjugate with respect to the left-right action of the double Weyl group 
  $W_L\times W_L$ on the set of roots  $\Delta(\ux^a)$.\\
 \Proof. 
\emph{ Item 1.}   Let  $u=1+X$ and $X\in \ux^a$. The associative algebra  $\ux^a$ is a subalgebra in the factor algebra of $\nx^a=\nx_+\ltimes \nx_-^a$ modulo the two sided ideal  
$\nx^a_-(L)$. There exists  $\overline{X}\in \nx^a$ such that  $X=\overline{X}\bmod \nx^a_-(L)$. It implies from \cite[Theorem 3.4]{P2} there exists  
$A,B\in \UT(n)\ltimes N_-^a$ such that $$A\overline{X}B=X_{\overline{D},\overline{\phi}}$$ for some rook placement  $\overline{D}\subset \Delta$ and a map  $\overline{\phi}:\overline{D}\to\Fq^*$. 
Decompose $A=A_1(1+\mu_L)$ and $B=(1+\nu_L)B_1$, where $A_1,B_1\in \UT(n)\ltimes U_-^a$ and $\mu_L, \nu_L\in \nx^a_-(L)$.  Since $\overline{X}=x+\la$ with $x\in \ux_+$ and $\la\in\ux_-^a$,  by Lemma \ref{muux}, we get  $(1+\mu_L)\overline{X}(1+\nu_L)=\overline{X}$. Then 
\begin{equation}\label{AXB}
A_1\overline{X}B_1=X_{\overline{D},\overline{\phi}}.
\end{equation} We denote $D=\overline{D}\cap \Delta(\ux^a)$ and $\phi =\overline{\phi}\vert_D$. Taking the equality  (\ref{AXB}) modulo the ideal  $\nx^a_-(L)$, we obtain $A_1XB_1=X_{D,\phi}$. Acting by a suitable  $h\in H$, we get $hX_{D,\phi}=X_D\in K(u)$, this proves 1).\\
\emph{ Item 2.} Suppose that  $AX_DB=X_{D'}$, where $A,B\in \GL^a(n)$. 
For any  $X\in \ux^a$ and  $1\leq k,m \leq n$,~ $k\ne m$,  we denote by   $r_{km}(X)$ rank of the  $I_k\times I_m$ block of the matrix $X$.  
By Lemma \ref{canon}, to prove statement 2), it is sufficient to show that  $r_{km}(X_D)=r_{km}(X_{D'})$ for all  $k\ne m$.
 We need the other realization of   $\gl^a(n)$ and also  $\ux^a$.

Consider the subalgebra  $\underline{\gl}^a(n)$ of $\mathrm{Mat}(2n)$ that consists of matrices  
$\left(\begin{array}{cc}A&B\\0&A\end{array}\right)$, where 
$A\in\px$. The algebra  $\gl^a(n)$ is isomorphic to the factor algebra of   $\underline{\gl}^a(n)$ modulo the ideal  $\left\{\left(\begin{array}{cc}0&B\\0&0\end{array}\right)\right\}$, where $B\in \px$.
Consider the decomposition of  segment  $[1,2n]$ into consecutive segments  $$I_1,\ldots,I_\ell, \tilde{I}_1,\ldots,\tilde{I}_\ell,$$ where $I_1,\ldots,I_\ell$ as above and for each  $1\leq m\leq \ell$  if $I_m=[a,b]$ then    $\tilde{I}_m=[a+n,b+n]$.
We realize the algebra  $\gl^a(n)$ as the subspace in  $\underline{\gl}^a(n)
$ consisting of matrices with zero  $I_k\times \tilde{I}_m$,~ $k\leq m$, blocks.  
\\
{\bf Example}. If  $n=2$ and   $\px$ is a Borel subalgebra, then
{\small{$$\gl^a(n) =\left\{\left(\begin{array}{cccc}
a_{11}&a_{12}&0&0\\
0&a_{22}&b_{22}&0\\
0&0&a_{11}&a_{12}\\
0&0&0&a_{22}\\ 
\end{array}\right)\right\},\quad \ux^a =\left\{\left(\begin{array}{cccc}
0&a_{12}&0&0\\
0&0&b_{22}&0\\
0&0&0&a_{12}\\
0&0&0&0\\ 
\end{array}\right)\right\}$$}}

Let $X\in \ux^a$. For any   $1\leq k<m\leq \ell $, we denote by  $R_{km}(X)$ the rank of  submatrix  of $X$ with the system of rows and columns $I_k\cup\ldots\cup I_m$.
For any  $1\leq m<k\leq \ell $, we denote by     $\tilde{R}_{km}(X)$ the rank of submatrix  of $X$ with the system of rows and columns $I_k\cup\ldots\cup I_\ell\cup\tilde{I}_1\cup\ldots\cup \tilde{I}_m $. Observe that  the ranks of  $R_{km}(X)$  and  $\tilde{R}_{km}(X)$ preserve if we replace  $X$ by  $AXB$ for any   $A,~B\in \GL^a(n)$.
We have $r_{km}(X_D)=r_{km}(X_{D'})$ for any  $k\ne m$. 
The canonical forms  $D$ and  $D'$ coincide (see Lemma \ref{canon}).  ~$\Box$

 The proof of Theorem  \ref{supclass} implies the following consequence.\\
 \Cor\Num. The elements  $u_D$ and  $u_{D'}$ belong to a common class if and only if  $r_{ij}(D)=r_{ij}(D')$ for all $ 1\leq i,~j\leq \ell$,~$i\ne j$.

Turn to a classification of the  $(\GL^a(n)\times \GL^a(n))$-orbits in $(\ux^a)^*$  and construction of supercharacters for $U^a$. 
We define the non-degenerate bilinear form on   $\ux^a$ by the formula 
\begin{equation}\label{form}
(X_1,X_2)=\la_2(x_1)+\la_1(x_2), 
\end{equation}
where $X_1=x_1+\la_1$ and $X_2=x_2+\la_2$.  One can directly calculate 
\begin{equation}\label{uniform}
(gX_1,X_2) = (X_1,X_2g),
\end{equation}
for all  $g\in G^a$. For any linear form  $\La\in(\ux^a)^*$, there exists $X_\La\in\ux^a$ such that 
$\La(Y)=(X_\La,Y)$ for any  $Y\in\ux^a$.  
  Taking into account formula (\ref{uniform}),  one can identify  the left-right  $(\GL^a(n)\times\GL^a(n))$-orbit of  $\La$ in  $(\ux^a)^*$  with the left-right orbit of  $X_\La$ in $\ux^a$.  

 For an arbitrary rook placement  $D$ and a map  $\phi: D\to\Fq^*$, we define the linear form
\begin{equation}\label{}
\La_{D,\phi} = \sum_{\gamma_+\in D_+}\phi(\gamma_+)E^*_{\gamma_+}  + 
\sum_{\gamma_-\in D_-}\phi(\gamma_-)F^*_{\gamma_-}.
\end{equation}
If $\phi$ is identically equal to one, we denote $\La_{D,\phi}=\La_D$.

Observe that  $\La_D(Y)=(X_{D^t},Y)$, where $D^t$ consists of all roots opposite to the roots from  $D$. It follows from Theorem \ref{supclass} that any   $(\GL^a(n)\times\GL^a(n))$-orbit in  $(\ux^a)^*$ contains a unique element of the form  $\La_D$. 

  Let us describe the stabilizer of  $\La_D$ with respect to the right action of  $U^a$ on $(\ux^a)^*$. 
 For each  $\gamma\in \Delta(\ux^a)$, consider the subset  $\Sc(\gamma)=\Sc_+(\gamma)\cup \Sc_-(\gamma)$, where $\Sc_\pm(\gamma)=\Sc(\gamma)\cap \Delta_\pm(U^a)$ is defined as follows:
 \begin{enumerate}
 	\item 
 	If  $\gamma\in I_k\times I_m$ for some   $k<m$ (i.e. $\gamma>0$), then  $\Sc_-(\gamma)=\varnothing$ and $\Sc_+(\gamma)$  consists of all roots $\al \in I_k\times I_t$ such that $\row(\al) =\row(\gamma)$ and   $k<t<m$. 
 	\item If  $\gamma\in I_k\times I_m$ for some   $k>m$ (i.e. $\gamma<0$), then $\Sc_-(\gamma)$ consists of all roots  $\al \in I_k\times I_t$ such that  $\row(\al) =\row(\gamma)$  and  $1\leq t < m$; $\Sc_+(\gamma)$ consists of all roots  $\al \in I_k\times I_t$ such that  $\row(\al) =\row(\gamma)$ and   $k< t \leq \ell$.
 	
 \end{enumerate}
 
 Let  $\ux^a_\gamma$ be the subspace spanned by the matrix units $E_\beta$,~  $\beta\notin \Sc_+(\gamma)$, and  $F_\beta$,~ $\beta\notin \Sc_-(\gamma)$. The subspace  $\ux^a_\gamma$ is a subalgebra  (more precisely, right ideal) in $\ux^a$.  Respectively,  $U^a_\gamma = 1+ \ux^a_\gamma$ is an algebra subgroup in $U^a$.
 
 Define $$\Sc(D)=\bigcup_{\gamma\in D} \Sc(\gamma), \quad\quad
 \ux^a_D= \bigcap_{\gamma\in D} \ux^a_\gamma, \quad\quad U^a_D =1+\ux^a_D= \bigcap_{\gamma\in D} U^a_\gamma.$$
 The next statement is proved by direct calculations.\\
 \Prop\Num. The subgroup $ U^a_D$ is the stabilizer of $\La_D$ with respect to the right action of  $U^a$ on $(\ux^a)^*$. 
 
  The formula  $$ \xi_D=\sum_{p\in L} \eps^{p\La_D(X)}, \quad u=1+X, \quad X\in \ux^a_D$$
 defines a character of the group  $U^a_D$. 
 Consider the character  $\chi_D$ of  $\GL^a(n)$ induced from the character  $\xi_{D}$  of the subgroup  $U^a_D$. 
 The following formula is true for   $\chi_D$:
\begin{equation}\label{chiUa}
\zeta_D(u) = \frac{|\La_D U^a|}{|\GL^a(n)\La_D\GL^a(n)|} \sum_{\mu\in \GL^a(n)\La_D\GL^a(n)}\eps^{\mu(X)}
\end{equation}
 where $u=1+X$  (see \cite[Corollary 3.1]{P1}).  Denote $K_D=1+\GL^a(n) X_D \GL^a(n)$.\\
 \Theorem\Num \label{ThSupUa}. The system of characters $\{\zeta_D\}$ and the decomposition of the group  $U^a$ into classes $\{K_D\}$ give rise to a supercharacter theory of the group  $U^a$. \\
 \Proof. It follows from Theorem  \ref{supclass} and \cite[theorem 3.3]{P1}. ~$\Box$

\subsection{ Supercharacter theory for  $\GL^a(n)$}\label{subsecgla}

In this subsection, following the general method of  paper \cite{P1}, we construct a supercharacter theory the group $\GL^a(n)$.

For each   $\gamma\in \Delta(\ux^a)$, we  consider the subgroup  $H_\gamma$ as follows. If $\gamma>0$ and $\gamma\in I_k\times I_m$, ~ $1\leq k<m\leq \ell$, then the subgroup $H_\gamma$ consists of all $h=(h_1,\ldots, h_\ell)\in L$,~  $h_i\in\GL(n_i)$ such that  $(h_k,\ldots,h_m)$ is a scalar submatrix.

If $\gamma<0$ and  $\gamma\in I_k\times I_m$, ~ $1\leq m<k\leq \ell$, then the subgroup   $H_\gamma$ consists of all  $h=(h_1,\ldots, h_\ell)\in L$ such that  $(h_1,\ldots,h_m,h_k,\ldots,h_\ell)$ is a scalar submatrix.\\
{\bf Example}. If $n=4$, ~ $P=B$,~ $\gamma_1=(1,3)$, and $\gamma_2=(3,1)$, then
{\small{$$H_{\gamma_1} =\left\{\left(\begin{array}{cccc}
		a&0&0&0\\
		0&a&0&0\\
		0&0&a&0\\
		0&0&0&b\\ 
		\end{array}\right)\right\},\quad H_{\gamma_2} =\left\{\left(\begin{array}{cccc}
		a&0&0&0\\
		0&b&0&0\\
		0&0&a&0\\
		0&0&0&a\\ 
		\end{array}\right)\right\}$$}}
For a rook placement  $D$, we define
$$H_D=\bigcap_{\gamma\in D}H_\gamma.$$
The next statement if proved by direct calculations.\\
\Prop\Num. For any rook placement  $D\subset \Delta(\ux^a)$, the subgroup $H_D$ consists of all  $h\in L$ such that $\Ad^*_h$ stabilizes all linear forms from $\GL^a(n)\La_D \GL^a(n)$. 
 
Consider the subset  $\Ac$ of all pairs  $(D, \theta)$, where $D$ is a rook placement contained in  $\Delta(\ux^a)$, and  $\theta$ is an irreducible character of the subgroup $H_D$.
Construct the subgroup  $G^a_D=H_DU^a_D$. 
The formula $$\xi_{\theta,\La_D}(hu)=\theta(h)\eps^{\La_D(X)}, \quad u=1+X, \quad X\in \ux^a_D$$
defines the character of  subgroup  $G^a_D$. 
Consider the character  $\chi_\al$ of  group  $\GL^a(n)$ induced from the character 
$$\sum_{p\in L} \xi_{\theta,p\La_D}$$
of  subgroup  $G^a_D$. 
There is a formula for  the character  $\chi_\al$  (see. \cite[Corollary 3.1]{P1}):
$$\chi_\al(g) = \frac{|L|^2}{|H_D|}\cdot\dot{\theta}(h)\zeta_D(u),$$
where $g=hu$, ~ $\zeta_D(u)$ are from  (\ref{chiUa}),~ $\dot{\theta}(h)=\theta(h) $ for $h\in H_D$ and equals to zero for  $h\notin H_D$. 

Turn to superclasses of  the group  $\GL^a(n)$.
Let  $h\in L$. Consider the smallest  $\GL^a(n)\times \GL^a(n)$-invariant ideal  $\ux^a_h$ in  $\ux^a$ such that $\Ad_h$ is identical on $\ux^a/\ux^a_h$. \\
{\bf Example}. Let $n=5$, ~$P=B$. We have
{\small{$$h =\left(\begin{array}{cccc}
		a&0&0&0\\
		0&a&0&0\\
		0&0&b&0\\
		0&0&0&a\\ 
		\end{array}\right), \mbox{where} ~a\ne b, \quad \ux^a_h =\left\{\left(\begin{array}{cccc}
		0&0&*&*\\
		\,*&0&*&*\\
		\,*&*&0&*\\
		0&0&*&0\\ 
		\end{array}\right)\right\}$$}}

 Consider the subset  $\Bc$ of all pairs $(D,h)$, where $D$ is a rook placement contained in  $\Delta(\ux^a)$, and  $h\in H_D$. Let $\omega_D$ be an orbit of the element $X_D$ with respect to  $\GL^a(n)\times \GL^a(n)$. 
We denote $$K_\beta= CL_L(h)\cdot (1+ \pi_h^{-1}\pi_h(\omega_D)),$$
where $CL_L(h)$ is the conjugacy class of  $h$ in $L$.\\
\Theorem\Num\label{GGsuper}. The system of characters  $\{\chi_\al:~ \al\in \Ac\}$ and the decomposition into classes  $\{K_\beta:~ \beta\in\Bc\}$ give rise to a supercharacter theory of the group  $\GL^a(n)$.\\
\Proof. It implies from  \cite[Theorem 3.3]{P1}. ~$\Box$

\section{Contractions of the orthogonal and symplectic groups}

As above  $\Fq$ is a finite field of $q$ elements and of characteristic  $\ne 2$. In this section, we use the notation  $\Gb=\GL(M, \Fq)$.
Respectively, for the subgroups in $\Gb$ we use new notations $\Pb$,~$\Lb$,~$\Hb$,~$\Ub_\pm$ and $\Nb_\pm$.
 The orthogonal (symplectic) subgroup  $G$  is defined by the involutive antiautomorphism   $X\to X^\dag$.

 Let  $M=2n+1$ for case  $B_n$  and $M=2n$ for $C_n$ and $D_n$. 
  Denote by $\Ib_M $ the $N\times N$ matrix in which all entries on the antidiagonal are equal to $1$ and all other entries are zeros. Consider the involutive antiautomorphism $X^\dag= \Ib_{M} X^t \Ib_{M}$ of the algebra of  $M\times M$ matrices over the field  $\Fq$. The orthogonal group $O(M)$ consists of all matrices  $g\in\Gb$ such that $g^\dag=g^{-1}$.   
  The symplectic group  $\Sp(M)$ is defined similarly applying the involutive automorphism   $X^\dag= \Jb_{2n}X^t\Jb_{2n}$, where $M=2n$ and $\Jb_{2n}=\left(\begin{array}{cc} 0&\Ib_n\\ -\Ib_n&0\end{array}\right).$    
  Respectively, the Lie algebra  $\gx$ of group  $G$ coincides with  the Lie algebra  $\{X\in\gl(M):~ X^\dag=-X\}$.

 We enumerate the rows and columns of  $M\times M$  matrices as follows:
 $$
 \begin{array}{c}
 n>\ldots > 1 > 0 > -1> \ldots >-n\quad \mbox{for} \quad M=2n+1,\\ n>\ldots > 1> -1> \ldots >-n \quad \mbox{for} \quad
 M=2n.
 \end{array}
 $$ 
 
Consider the partition  $\Ic$ of  segment   $I=[-n,n]$ in case  $B_n$ (respectively, $I=[-n,n]\setminus \{0\}$ in cases $C_n$ and $D_n$) into  the system of consecutive segments.   We enumerate the  components of partition  $\Ic=\{I_{\ell},\ldots, I_0,\ldots , I_{-\ell}\}$ comparably with the enumeration of rows and columns. Suppose that the partition $\Ic$  is symmetric relative to zero. By the partition $\Ic$, we construct the parabolic subgroup $\Pb=\Lb\Ub_+$ and the complementary unipotent subgroup  $\Ub_-$. As above  $\Ub_\pm =1+\Uc_\pm$ are the algebra subgroups.

Intersection of these subgroups with $G$ defined the parabolic subgroup  $P=LU_+$ and the complementary unipotent subgroup $U_-$. 
The subgroup $L$ consists of all matrices  $\mathrm{diag}(A_{\ell},\ldots, A_0, \ldots, A_{-\ell})$, where for each  $\ell\geq i\geq {-\ell}$ the transpose of  matrix   $A_i$ about its anti-diagonal coincides with  $A_{-i}^{-1}$. 
We preserve notations  $\ux_\pm$ for Lie algebras of unipotent groups  $U_\pm$.
 We identify  $\ux_-$ with  $\ux^*$.

We extend the involutive antiautomorphism of  subgroup $\Pb$  to the one of $\Gb^a=\Pb\ltimes \Ub_-^a$ setting  $(1+\la)^\dag=1-\la$, ~ $\la\in\Uc_-^a$, and $ (gv)^\dag=v^\dag g^\dag$, $g\in\Pb$, ~$v=1+\la$. The parabolic contraction  $G^a$ of orthogonal and symplectic groups consists of all elements  $gv\in\Gb^a$, ~ $g\in \Pb$ and  $v\in \Ub_-^a$ obeying  $(gv)^\dag=(gv)^{-1}$.  

 The group  $G^a$ is a semidirect product of  $L$ and the unipotent subgroup  $U^a=U\ltimes U_-^a$ with Lie algebra  $\ux^a = \ux_+\oplus \ux_-^a$. 

The group  $\Gb^a$   acts on $\ux^a$ by the formula 
\begin{equation}\label{gbaction}
 g\centerdot X= gXg^\dag,
\end{equation} 
where $g\in \Gb^a $ and  $X\in\ux^a$. We define the action of    $\Gb^a $ on  $(\ux^a)^*$ by  $$g\centerdot \la(X)=\la(g^\dag X g).$$

Consider the Cayley map   $f:\Ub^a\to\Uc^a$; its value  $f(u)$ at  $u=1+X\in\Ub^a$ is calculated by the formula 
 $$f(u) = \frac{2(u-1)}{u+1} = \frac{2X}{X+2} = \sum_{k=0}^\infty(-1)^k\frac{X^{k+1}}{2^k}.$$

\subsection{Supercharacter theory for orthogonal and symplectic  $U^a$}

In this section, we classify $\Gb^a$-orbits in $\ux^a$ and $(\ux^a)^*$ and  construct supercharacter theory for  $U^a$. We need definition of 
rook placements  for orthogonal and symplectic groups.

Let  $\Delta$ be the root system of the orthogonal or symplectic group  $G$.  Ordering of the segment  $I=[-n,n]$ defines the system  $\Delta_+$  of positive roots  (respectively, the system of negative roots  $\Delta_-$).

With each root $\al$, we associate some element of $I\times I$ according to the following rule.
For positive roots  $\al=\epsilon_i-\epsilon_j $, ~$i>j$,  and  $\al=\epsilon_i+\epsilon_j $  we associate the pairs   $(i,j)$ and $(i,-j)$ respectively.
 Analogously, for   $\al=\epsilon_i$,  it is the pair $(i,0)$,  for  $ \al=2\epsilon_i$ -- the pair  $(i,-i)$.
With the negative root $-\al<0$ we associate the pair transpose to the pair of  $\al$.
 For example, for   $\al=-\epsilon_i-\epsilon_j $,  it is    $(-j,i)$.
 For simplicity,   if  $\al=\epsilon_i-\epsilon_j $, ~$i>j$,  we write  $\al=(i,j)$ (similarly in other cases).

 For each root $\al=(i,j)$, we  consider the pair  $\al'=(-j,-i)$ (observe that the last pair does not associate with any root,  except $\al=\al'=(i,-i)$).
 
To each root $\al\in \Delta_+$, we attach the root vector  $\Ec_\al=E_\al+\eta(\al)E_{\al'}$, where $\eta(\al)\in\{1,0,-1\}$.
Respectively, to a root  $\al\in \Delta_-$, we attach the root vector  $\Fc_\al=F_\al+\eta(\al)F_{\al'}$. The system of root vectors  $\{\Ec_\al\}\cup\{\Fc_\al\}$  is a basis in $\ux^a$.
The relative elements of the dual basis we denote by 
$\Ec^*_\al$ and  $\Fc^*_\al$.  

Similarly to the previous section, we define $\Delta_+(\ux^a)$ that consists of all  $\al\in\Delta_+$ such that  $\Ec_\al\in\ux_+$. Respectively, $\Delta_-(\ux^a)=\{ \al\in\Delta_-:~ \Fc_\al\in\ux^a_-\}$. We define $\Delta(\ux^a)=\Delta_+(\ux^a)\cup \Delta_-(\ux^a).$

 \Def\Num. We say that a subset   $D\subset \Delta$ is a rook placement if there is at most one pair  $(i,j)\in D\cup D'$   in any row  and column. 
 In other words, $D$ is a rook placement in orthogonal or symplectic case if and only if  
 $D\cup D'$ is a rock placement in the root system of  $\GL(M)$.

\Ex~\Num\label{exone}. Let $G=C_2$, ~$\Delta_+=\{\al_1,\al_2,\al_1+\al_2,~2\al_1+\al_2\}$,~ $D=\{\al_1+\al_2, -\al_2\}$. On the diagram the pairs from  $D\cup D'$ are marked by $\otimes$. On the diagonal we put the numbers of relative rows  (columns).

{\small $$\left(\begin{array}{cccc}
2&\cdot&\otimes&\cdot\\
\cdot&1&\cdot&\otimes\\
\cdot&\otimes&-1&\cdot\\
\cdot&\cdot&\cdot&-2
\end{array}\right)$$.}

A rook placement  $D\subset \Delta$ splits  $D=D_+\cup D_-$, where $D_\pm=D\cap\Delta_\pm$.

 To the rook placement $D\subset\Delta$ and   map  $\phi:D\to\Fq^*$, we attach the element  $\Xc_{D,\phi}\in\ux^a$ and  $\La_{D,\phi}\in(\ux^a)^*$ defined by the formulas
$$ \Xc_{D,\phi}=\sum_{\gamma\in D_+}\phi(\gamma)\Ec_\gamma + \sum_{\gamma\in D_-}\phi(\gamma)\Fc_\gamma,$$
$$ \La_{D,\phi}=\sum_{\gamma\in D_+}\phi(\gamma)\Ec^*_\gamma + \sum_{\gamma\in D_-}\phi(\gamma)\Fc^*_\gamma.$$ 

Fix the non-square element $\delta\in \Fq^*$. 
In the case  $C_n$, we say that a root  $\gamma$ belongs to anti-diagonal if $\gamma=(i,-i)$ for some  $n\geq i\geq -n$.\\
In the case  $B_n$ and  $D_n$, we refer to pair $(D,\phi)$ as a \textit{ basic pair},
if $\phi(\gamma)=1$ for any  $\gamma\in D$.\\
\Def\Num\label{special}. In the case $C_n$, we refer to pair  $(D,\phi)$ as a basic pair if   \\
1)~ $\phi(\gamma)\in\{1,\delta\}$. \\
2) ~  $\phi(\gamma)=\delta$ implies   $\gamma$ belongs to anti-diagonal.\\
3) ~ For any $\ell\geq k \geq -\ell$, the subset  $D\cap (I_k\times I_{-k})$ is contained in the anti-diagonal, and
 there is at most one root  $\gamma$ in this subset such that   $\phi(\gamma)=\delta$.

By definition,  we let  $d_k(\phi)$  equal to $-1$  if there is a root from 3), otherwise we take $d_k(\phi)=1$. 
Below we introduce  notations that are similar to the relative notations for series $A$ (the difference in substitution $D\cup D'$ for $D$  and in changing the row and column enumeration). For any rock placement  $D\subset \Delta$ and pair of segments  $I_k, I_m$, where $k \ne m$ and $k+m\geq 0$, we define 
$$r_{ij} = r_{ij}(D) =  |(D\cup D')\cap(I_k\times I_m)|,$$
$$n'_{km} = n-\sum_{s>k}n_s -\sum_{t>m}r_{kt},$$
$$n''_{km} = n-\sum_{s>m}n_s -\sum_{t>k}r_{tm}.$$
\Def\Num.  We say that a basic pair  $(D_c,\phi_c)$  has canonical type if for each pair ~ $k\ne m$, ~$k+m \geq 0$, the subset  
$D_c\cap(I_k\times I_m)$ coincides with the subset 
\begin{equation}\label{knem}
\{ (n_{km}'-1, n_{km}''-1),\ldots, (n_{km}'-r_{km}, n_{km}''-r_{km})\},
\end{equation}
in the cases  $B_n$, $D_n$ and  $C_n$ (for $k+m>0$). In the case  $C_n$ and $m=-k$, the subset  $D_c\cap(I_k\times I_{-k})$  coincides with  
\begin{equation}\label{kminusk}
\{(n_{k,-k}'-1,-n_{k,-k}'+1)),\ldots, (n_{k,-k}'-r_{k,-k},-n_{k,-k}'+ r_{k,-k})\},
\end{equation} 
furthermore the values of $\phi_c$ at all roots from  (\ref{kminusk}) apart from possibly the first root  equal to one. 

Let $W_\Lb$ be the Weyl group of the subgroup  $\Lb$ in $\Gb$.  The action of  $W_\Lb$ on  $\ux^a$ is obtained by restriction of the $\Gb^a$-action  on $\ux^a$ (see formula (\ref{gbaction})).  Under this action the elements  $\Xc_{D,\phi}$ transform to elements of the same type;  this fact defines an action of  $W_\Lb$ on the set of basic pairs $\{(D,\phi)\}$. 

\Lemma\Num\label{canontwo}. For any basic pair $(D,\phi)$, there exists an element  $w\in W_\Lb$ such that  $w\centerdot (D,\phi)$ has the canonical type  $(D_c,\phi_c)$.  The canonical  pair $(D_c,\phi_c)$ is uniquely determined by   $(D,\phi)$.\\
\Proof. We proceed the proof for $D_n$.  The other series are treated similarly. We prove by the induction method on  $n$.

Let $\Dc=D\cup D'$. Decompose the square  $I\times I$ into four parts 
$I\times I=I_{11}\cup I_{12}\cup I_{21}\cup I_{22}$, where
$I_{11} =[0,n]\times [1,n]$,~ $I_{12} =[0,n]\times [-n,0]$,~
$I_{21} =[-n,-1]\times [1,n]$,~ $I_{22} =[-n,-1]\times [-n,0]$.

We take $\Dc_{ij}=\Dc\cap I_{ij}$. Observe that  
$\Dc_{11}'= \Dc_{22}$,~ $\Dc_{12}'=\Dc_{12}$ and $\Dc_{21}'=\Dc_{21}$.

 Each element $w\in\Lb$ has the form  $w=\mathrm{diag}(A_{\ell},\ldots, A_0, \ldots, A_{-\ell})$, where $A_i\in\GL(n_i)$. Decompose $w=w_1w_2$, 
 where $w_1= \mathrm{diag}(A_{\ell},\ldots, A_0, E,\ldots, E)$ and 
 $w_2=\mathrm{diag}(E,\ldots, E, A_{-1}, \ldots, A_{-\ell})$.
 
 Then $w\centerdot \Dc_{11}=w_1\Dc_{11}w_2^\dag$, ~ $w\centerdot \Dc_{12}=w_1\Dc_{12}w_1^\dag$, ~ $w\centerdot \Dc_{21}=w_2\Dc_{21}w_2^\dag$.
Replacing  $\Dc$ by suitable    $w\centerdot\Dc$,  we consider  $\Dc_{11}$  to be of canonical type in $I_{11}$ (see Lemma \ref{canon}).

Suppose that  $\tilde{I}_{12}$ is a square obtained by deleting in  $I_{12}$ the rows from 
$\row(\Dc_{11})$ and columns from  $\col(\Dc_{22})$ (here the number of rows equals to the number of columns as $\Dc_{22}=\Dc_{11}'$).  
Since there are no roots from  $\Dc_{12}$ in rows  from  $\row(\Dc_{11})$ and columns from  $\col(\Dc_{22})$, we have $\Dc_{12}\subset  \tilde{I}_{12}$. 
Applying the induction assumpsion, we can replace  $\Dc_{12}$ by some $w_1\Dc_{12}w_1^\dag$ of canonical type, where 
$w_1$ stabilizes the rows from  $\row(\Dc_{11})$.
Analogously for  $w_2$ and $\Dc_{21}$.

 Uniqueness follows from the fact that  the numbers  $\{r_{km}\}$ does not change under  the action of   $\Lb$. ~$\Box$
 
For parabolic subgroups, the next theorem is proved in   \cite[Proposition 3.3]{P6}.\\
\Theorem\Num\label{classUaos}. 1)  Each  $\Gb$-orbit in  $\ux^a$ contains an element  of the form   $ \Xc_{D,\phi}$  for some  basic pair $\{D,\phi\}$.\\
2)~Two elements $\Xc_{D,\phi}$ and $  \Xc_{D',\phi'}$ for basic pairs 
$(D,\phi)$ and $(D',\phi')$ belong to a common $\Gb^a$-orbit if and only if   
$(D,\phi)$ and $(D',\phi')$ belong to a common  $W_\Lb$-orbit.\\
\Proof.  \emph{Item 1}. Our goal is to prove statement  1) for the Borel contraction. In this case $P=B=HN$ is a Borel subgroup in the orthogonal or symplectic group, and $\Pb=\Bb=\Hb \Nb$ is a Borel subgroup in  $\Gb$.
 The unitriangular  group  $\Nb$ is an algebra subgroup  $\Nb=1+\Nc$.
The group  $\Gb^a$ acts on  $\Nc^a=\Nc_+\ltimes \Nc_-^a$.    
   To prove statement 1), we use the method analogous to \cite[Proposition  3.3]{P2}. \\ 
1) Let  $X=x+\la\in \nx^a$,  where $x\in \nx_+$ and  $\la\in\nx_-^a$.
 We consider  $X,~x,~\la$ as the elements from  $\Nc^a$  obeying $Y^\dag=-Y$.  
  It follows from   \cite{AN-1,AN-2} that there exists   $a\in \Nb_+ $ such that $axa^\dag=\Xc_{D_+,\phi_+}$, where  $D_+$  is a rock placement in $\Delta_+$ and $\phi:D_+\to \Fq^*$. 
Take $x_+=\Xc_{D_+,\phi_+}$. We obtain  $aXa^\dag=x_++a\la a^\dag$.\\
2) In what follows  $X=x_++\la$.  In this subitem, we consider decomposition   (\ref{nminusa}) for $\Nc_-^a=\Nc_+^*$   and dual decomposition  $\nx_-^a$ (see (\ref{fspace}) and (\ref{sspace})); we show that  applying $\Nb^a$-action we get a similar element with  $\la$  in the subspace  (\ref{sspace}). As above $\Dc=D\cup D'$.

Consider  decomposition $\Nc_+=\Lc_+\oplus \Lc_0$ into direct sum of two left ideals
$$\Lc_+= \spann \{E_{ij}:~ j\notin \row(\Dc_+)\}=\{y\in\Nc_+:~ yx_+=0\},$$
$$
\Lc_0=\spann \{E_{ij}: ~ j\in \row(\Dc_+)\}.$$
Then $\Nc_-^a=\Lc_+^\perp\oplus\Lc_0^\perp$, where 
$$
\Lc_+^\perp=\spann\{F_{km}:~ k\in \row(\Dc_+)\}=x_+\Nc_-^a, \quad \Lc_0^\perp=\spann\{F_{km}:~ k\notin \row(\Dc_+)\}
$$
We have $\Nc_-^a=x_+\Nc_-^a\oplus\Lc_0^\perp$.

Similarly, consider decomposition  $\Nc_+=\Rc_+\oplus \Rc_0$ into  direct sum of two right ideals, where $\Rc_+=\Lc_+^\dag$ and $\Rc_0=\Lc_0^\dag$.

Since the subspaces  $\Lc_+$, $\Lc_0$, $\Rc_+$, and $\Rc_0$ are spanned by vectors from standard basis  $\{E_{ij}\}$ in  $\Nc_+$,  we obtain 
\begin{equation}\label{nminusa}
\Nc_-^a=\left(x_+\Nc_-^a + \Nc_-^ax_+\right)\oplus \left(\Lc_0^\perp\cap \Rc_0^\perp\right)
\end{equation}

Both subspaces in this sum are invariant under action of involutive anti\-auto\-morphism. Therefore, 
$\nx_-^a$ is a direct sum of intersections of these subspaces with  $\nx_-^a$: 
\begin{equation}\label{fspace}
 \left(x_+\Nc_-^a + \Nc_-^ax_+\right)\cap\nx_-^a= \left\{x_+Q^\dag+Q x_+:~ Q\in \Nc_-^a \right\},
\end{equation}
\begin{equation}\label{sspace} \left(\Lc_0^\perp\cap \Rc_0^\perp\right)\cap \nx_-^a=
\left\{V\in \nx_-^a:~ V(\Lc_0)=0\right\}=\left\{V\in \nx_-^a:~ V(\Rc_0)=0\right\}.\end{equation}

Return to the element    $X=x_++\la\in \nx$. The linear form  $\la\in\nx_-^a$  splits  
$\la=x_+Q^\dag +Q x_+ + \la_0$, where $\la_0$ belongs to subspace  (\ref{sspace})  in decomposition $\nx_-^a$.
We get 
$$(1-Q)X(1-Q^\dag) =  x_+ - x_+Q^\dag -Qx_+ +\la= x_++\la_0.$$
3) So, we can assume that  $X=x_++\la_0$, where $\la_0(\Lc_0)=\la_0(\Rc_0)=0$ (recall that  $\Rc_0=\Lc_0^\dag$).  Let us show that in $\Gb^a
$-orbit of the element $X$ there is an element of the form   $x_+ +\la'$, where $$\la'=\Xc_{D_-',\phi_-'}$$ for some rook placement   $D_-'\subset \Delta_-$ and $\phi_-': D_-' \to \Fq^*$.

For $\la_0$  there exists  $a\in \Nb_+$ such that  $a\la_0a^\dag=\la'$, where $\la'$ has the mentioned above form. Let us show that we can chose  $a\in \Nb_+$ such that  $ax_+a^\dag=x_+$.

Consider the algebra subgroups $L_+=1+\Lc_+$ and  $L_0=1+\Lc_0$ in $\Nb_+$.
Since $\Nc_+$ is a direct sum of left ideal  $\Nc_+=\Lc_++\Lc_0$, each  $a\in \Nb_+$ can be uniquely presented in the form  $a=a_+a_0$, where $a_+\in L_+$ and $a_0\in L_0$.

By definition of $\Lc_+$, we have $a_+x_+=x_+$. From $x_+^\dag=-x_+$, we get $x_+a_+^\dag=x_+$.

Let us show that $a_0\la_0=\la_0$ for  any  $a_0\in L_0$. Indeed,   $a_0=1+y_0$, where $y_0\in\Lc_0$. Then $a_0\la(x)=
\la_0(xa)=\la_0(x)+\la_0(xy_0)$. Since $\Lc_0$ is a left ideal  in  $\Nc_+$, it follows $\la_0(xy_0)=0$. This proves  $a_0\la_0=\la_0$.
From  $\la_0^\dag=-\la_0$, we obtain $\la_0a_0^\dag=\la_0$. We have  $$a_+Xa_+^\dag=a_+x_+a_+^\dag + a_+\la_0a_+^\dag = x_+ + a_+a_0\la a_0^\dag a_+^\dag = x_++a\la_0a^\dag = x_++\la'.$$

4)  So, each  $\Nb^a$-orbit of the element  $X\in \nx^a$ has an element  $X'=x_++\la'$, where $x_+=\Xc_{D_+,\phi_+}$ and $\la'=\Xc_{D_-',\phi_-'}$.  Finally, as in subitem  3), replacing  $X'$ by a suitable element of the form  $(1+Q)X'(1+Q^\dag)$, we can delete rocks of  $\Dc_-'$ belonging to rows and columns from  $\Dc_+$.
We obtain the element   $\Xc_{D,\phi}$, where $D=(D_+,D_-)$ is a rook placement in  $\Delta$ and  $\phi$ is a map  $D\to\Fq^*$ coinciding with $\phi_+$ on $D_+$  and with $\phi_-'$ on $D_-$.
Finally, replacing  $X\to hXh^\dag$ by a suitable  $h\in \Hb$, we make $(D,\phi)$ to be a basic pair.
The statement if Item 1 is proved. \\
\textit{Item 2.}  The proof of statement 1) for parabolic contraction is proceed similarly to Item 1 from Theorem
 \ref{supclass}.
One can prove  statement  2)  similarly to Item 2 of Theorem  \ref{supclass}  applying Lemma  \ref{canontwo}. 
 ~$\Box$ \\
\Cor\Num.
Two elements  $\Xc_{D,\phi}$ and $  \Xc_{D',\phi'}$ for basic pairs 
$(D,\phi)$ and $(D',\phi')$ belong to a common  $\Gb^a$-orbit if and only if   $r_{km}(D)=r_{km}(D)$
for all  $\ell\geq k, m\geq -\ell$, ~$k\ne m$,  and $d_k(\phi)=d_k(\phi')$ for all  $\ell\geq k \geq -\ell$,~ $k\ne 0$.\\
\Proof.  Analogously to Theorem  \ref{supclass} and  \cite[Proposition 3.3]{P6}.
~$\Box$ 

Using identification of  $\ux^a$ and $(\ux^a)^*$, we obtain the following statement.\\
\Prop\Num\label{dualclassUaos}. 1)  Each  $\Gb$-orbit in  $(\ux^a)^*$ contains an element of the form   $ \La_{D,\phi}$  for some basic pair  $\{D,\phi\}$.\\
2) Two elements  $\La_{D,\phi}$ and $  \La_{D',\phi'}$ for basic pairs 
$(D,\phi)$ and  $(D',\phi')$ belong to a common  $\Gb^a$-orbit if and only if   
$(D,\phi)$ and  $(D',\phi')$ belong to a common  $W_\Lb$-orbit.\\
 \Def\Num. We say that element  $u$ and $u'$ from $U^a$ belong to a common class if
 the elements  $f(X)$ and $f(X')$ belong to a connon $\Gb^a$-orbit in  $\Uc^a$ ( here  $f$ is the Cayley map).
 
 Denote $u_{D,\phi}=f^{-1}(\Xc_{D,\phi})$. 
 Theorem  \ref{classUaos} yields a classification of classes in the group  $U^a$.
 
 For each basic pair  $(D,\phi)$ we consider a function on $U^a$ of the form
  \begin{equation}\label{charUa}
 \sigma_{D,\phi}(u) =  \sum_{\mu\in \Gb^a\centerdot \La_{D,\phi}}\eps^{\mu(f(u))}
 \end{equation}
    Denote by $K_{D,\phi}$ the preimage of $\Gb^a$-orbit of  $\Xc_{D,\phi}$ under the Cayley map.\\ 
   \Theorem\Num \label{SupUaos}. The system of functions  $\{\sigma_{D,\phi}\}$ and the decomposition of group  $U^a$ into classes  $\{K_{D,\phi}\}$ give rise to a supercharacter theory of the group  $U^a$. \\
 \Proof.  It follows from   \cite{AFN} that  the system of functions 
 $$ \sigma_{\la}(u) =  \sum_{\mu\in \Ub^a\centerdot \la}\eps^{\mu(f(u))}$$
 and  subsets $K_x= f^{-1}(\Ub^a\centerdot x)$, where $\la$ and  $x$ run through  the  systems of representatives of  $\Ub^a$-orbits in  $(\ux^a)^*$ and  $\ux$ respectively, give rise to a supercharacter theory of the group  $U^a$. In particular, the functions   $\sigma_{\la}$ are characters of some representations  and pairwise orthogonal (disjoint). Therefore,  $ \sigma_{D,\phi}$ are also characters and pairwise orthogonal. Theorem   \ref{classUaos} implies that the decomposition into subsets   $K_{D,\phi}$   form a partition of the group  $U^a$ more coarse than the decomposi\-tion  $\{K_x\}$. Finally, if  $u$ and  $u'$ belong to a common class, then there exists  $g\in \Gb^a$, such that  $f(u)=g\centerdot f(u')$. Substituting for (\ref{charUa}), we verify that the values of  $\sigma_{D,\phi}$ at $u$ and $u'$ coincide.  ~$\Box$

\subsection{Supercharacter  theory for orthogonal and symplectic  $G^a$}

In subsection  \ref{subsecgla}, by a pair  $\gamma=(i,j)$, ~$i\ne j$, we  construct the subgroup  $H_\gamma$ in $L\subset \GL(n)$.  In notations of this section, it is the subgroup  $\Hb_\gamma$  in the block-diagonal group $\Lb$.

Let  $D$ be a rook placement in  $\Delta(\ux^a)$. For each root  $\gamma$, we define  the subgroup  $H_\gamma$ in $L$ as  intersection of  $\Hb_\gamma\cap\Hb_{\gamma'}$ with orthogonal or symplectic group $G$.

Define  $$H_D=\bigcap_{\gamma\in D}H_\gamma.$$

Denote by $\ux_D^a$ the intersection of all  $\mathrm{Ker}\, g\centerdot\La_{D,\phi}$, where $g\in\Gb^a$. The subset  $\ux_D^a$ is a  $\Gb_a$-invariant ideal in  $\ux^a$.\\
\Lemma\Num\label{hddd}. For each  $h\in H_D$ the operator  $\Ad_h$ is identical on $\ux^a/\ux^a_D$.\\
\Proof. Consider the subspace  $W_D^*\subset (\ux^a)^*$ spanned by  $g\centerdot\La_{D,\phi}$, where $g\in\Gb^a$.  The subspace  $W_D^*$ is invariant with respect to action of  $\Lb$. Therefore,  $W_D^*$ decomposes into  $\Lb$-irreducible components  $W_{km}^*$ related to segments of decomposition of  $I=[-n,n]$.  If $\gamma \in I_k\times I_m$,~ $k>m$, then  $W_D^*$  contains all components  $W_{ij}^*$,~ $k\geq i>j\geq m$.  
 If $\gamma \in I_k\times I_m$,~ $k<m$, then   $W_D^*$  contains all components   $W_{ij}^*$,~ $i\leq k <m\leq j$.  For each  $h\in H_\gamma$, the operator  $Ad^*_h$  is identical on each mentioned above  compotent $W_{ij}^*$.
The subspace $W_D^*$ is a direct sum of all enumerated components  $W_{ij}^*$ constructed for  roots $\gamma\in D$.
Therefore, for each  $h\in H_D$ the operator $Ad^*_h$ is identical on   $W_D^*$. It implies that for  $h\in H_D$ the operator  $Ad_h$ is identical on   $\ux^a/\ux^a_D$. ~$\Box$\\ 
\Ex\Num. For root system  $D=\{\al_1+\al_2, -\al_2\}$ of $G=C_2$ (see example \ref{exone}), we have  $W^*_D=\mathrm{span}\{\Fc^*_{-\al_2},~\Fc^*_{-\al_1-\al_2}, ~ \Fc^*_{-2\al_1-\al_2},~ \Ec^*_{\al_1},~ \Ec^*_{\al_1+\al_2}, \Ec^*_{\al_2}\},$
$$\ux_D^a = \mathrm{span}\{\Fc_{-\al_1},~  \Ec_{2\al_1 + \al_2}\},\quad \quad H_D = \{\pm E\}.$$
On this diagram, the roots of   $D\cup D'$ is marked as  $\otimes$. The squares relative to basic elements of  $\ux_D^a$ are marked by  $\bullet$.
On diagonal we put the numbers of relative rows (columns).  

{\small $$\left(\begin{array}{cccc}
2&\cdot&\otimes&\bullet\\
\bullet&1&\cdot&\otimes\\
\cdot&\otimes&-1&\cdot\\
\cdot&\cdot&\bullet&-2
\end{array}\right)$$.}

 Consider the set  $\Ax$ of all triples  $\ax=(D,\phi,\theta)$, where $(D,\phi)$ is a basic pair from $\Delta(\ux^a)$,  and $\theta$ in an irreducible representation of the subgroup $H_D$. Below we attach to each  $\ax$ the character  $\sigma_\ax$ of subgroup $G^a$.
  
  Denote by  $U^a_D$ the normal subgroup  $f^{-1}(\ux^a_D)$ in $U^a$. Lemma  \ref{hddd} implies that the subgroup  $G_D^a = H_D\ltimes U_D^a$ is abelian modulo    $U^a_D$.  
  The characters  $\sigma_{D,\phi}$ defined by formula  (\ref{charUa})  can be taken modulo  $U^a_D$ (see formula  (\ref{charUa})). The formula $$\xi_\ax (g) = \theta(h) \sigma_{D,\phi}(u),$$
  ~ $g=hu$,~$h\in H_D$,~ $u\in U^a$, defines the character of subgroup  $G_D^a$.
  Consider the character  
  $$\sigma_\ax = \Ind(\xi_\ax, G_D^a, G^a).$$ 
 
 We turn to construction of superclasses. 
For each  $h\in L$,
we consider the smallest  $\Gb^a$-invariant ideal  $\ux^a_h$ in  $\ux^a$ such that   $\Ad_h$ is identical on  $\ux^a/\ux^a_h$. Let  $U_h^a$ be  the corresponding normal subgroup in $U^a$, and     $\Pi_h$ be the natural projection  $U^a$ onto $U^a/U^a_h$.

Consider the set $\Bx$ of all triples  $\bx=(D,\phi,h)$, where  $(D,\phi)$ is a basic pair in  $\Delta(\ux^a)$  and  $h\in H_D$. 
By each triple $\bx$, we define the class $$K_\bx= \mathrm{CL}_L(h)\cdot  \Pi_h^{-1}\Pi_h \left(K_{D,\phi}\right),$$
where $\mathrm{CL}_L(h)$ is the conjugacy class of  $h$ in $L$, and  $K_{D,\phi}$ from Theorem \ref{SupUaos}.\\
\Theorem\Num\label{Gsuper}. The system of characters   $\{\sigma_\ax:~ \ax\in \Ax\}$ and decomposition of  $G^a$ into  classes $\{K_\bx:~ \bx\in\Bx\}$ give rise to a supercharacter theory of the orthogonal or symplectic   $G^a$.\\
\Proof.  The characters  $\{\sigma_\ax\}$ are orthogonal because  
$$
\sigma_\ax(g) = c\, \dot{\theta} (h) \sigma_{D,\phi}(u),
$$
where $\sigma_{D,\phi} $ defined by formula (\ref{charUa}), ~$g=hu$,~ $h\in L$,~ $u\in U^a$, and  $\dot{\theta}(h)=\theta(h) $ for $h\in H_D$ and is zero for $h\notin H_D$. 

Let us show that  $\sigma_\ax$ is constant on  $K_\bx$, where
$\bx=(D',\phi',h)$, ~$h\in H_{D'}$.
If $h\notin H_D$, then $\sigma(K_\bx) = 0$.

Suppose $h\in H_D$. Then $\ux^a_h\subset \ux_D^a$ and $U_h^a\subset U_D^a$. 
The  $\sigma_\ax$ and class  $K_\bx$ can be taken modulo  $U^a_h$.
The character  $\sigma_\ax$ is constant on class $K_\bx$ because 
$\sigma_{D,\phi}$ is constant on  $K_{D',\phi'}$ for $U^a/U^a_h$.~$\Box$

\section{Funding}
The paper is supported by the RFBR grant 20-01-00091a.


\begin{thebibliography}{99}
	


\bibitem{DI}
Diaconis P., Isaacs I.M., \textit{Supercharacters and  superclasses for
	algebra groups.},Trans.Amer.Math.Soc., {\bf 360} (2008), 2359-2392.


\bibitem{A1}
Andr\'{e} C.A.M., \textit{Basic characters of the unitriangular group}, J. Algebra, {\bf 175} (1995), 287-319.
\bibitem{A2}
Andr\'{e} C.A.M., \textit{Basic sums of coadjoint orbits of the unitriangular group},  J. Algebra,
{\bf 176}(1995), 959-1000.
\bibitem{A3}
Andr\'{e} C.A.M., \textit{The basic character table of the unitriangular group},  J. Algebra,
{\bf 241} (2001), 437-471.

\bibitem{AN-1}
Andr\'{e} C.A.M.,  Neto  A.M., \textit{Supercharacters of the finite the Sylow subgroup of the finite symplectic and orthogonal groups},  Pacific Math. Journal, {\bf  239} (2009), 201-230.
\bibitem{AN-2}
Andr\'{e} C.A.M.,  Neto  A.M., \textit{A supercharacter theory for the Sylow p-subgroups of the finite symplectic and orthogonal groups},  J. Algebra, {\bf  322} (2009), 1273-1294.
\bibitem{AFN} Andr\'{e} C.A.M.,  Freitas J.P.,   Neto A.M.,
\textit{A supercharacter theory for involutive algebra groups}, J. Algebra,
{\bf 430} (2015), 159-190.


\bibitem{IW} In\"{o}n\"{u} E., Wigner E.P., \textit{On the contraction of groups and their representations}, Proc. Nat-Acad. Sci., {\textbf 39} (1953), 510-524. 


	\bibitem{GOV}
Gorbatsevich 	V.V.,   Onishchik A.L., Vinberg  E.B. (1994).\textit{ Lie groups and Lie algebras III.} Encyclopaedia Math. Sci., vol. 41, Berlin: Springer.



\bibitem{P1}
Panov A.N.   Supercharacter theories for algebra group extensions. J. Group theory, {\bf 24} (2021), 263-283.

\bibitem{P2} Panov A.N., \textit{Supercharacter theory for  the Borel group contraction of  $\GL(n,\Fq)$}, Vestnik St. Peterburg Univ. Math.,  \textbf{53} (2020),  no 2, 162-173. 

\bibitem{P3} 
Panov A.N. (2016).  \textit{Supercharacter theory for groups of invertible elements of reduced algebras}, St. Petersburg Math. Journal,    \textbf{27} (2016), 1035-1047.

\bibitem{P4}
Panov A.N., \textit{Supercharacters for  finite groups of  triangular type}, Comm.Algebra, {\bf 46} (2018), no.3, 1032-1046.

\bibitem{P5}
Panov A.N., \textit{Towards a supercharacter theory of  parabolic subgroups
	in $\GL(n,\Fq)$}, Algebr. Represent. theory, {\bf 21} 2018,  no. 5,  1133–1149.

\bibitem{P6}
Panov A.N., \textit{Two supercharacter theories  for the parabolic subgroups in orthogonal and symplectic groups},  J. Algebra,  \textbf{539} (2019),  37-53. 

\bibitem{Thiem}
Thiem N., \textit{Supercharacter theories of type $A$ unipotent  radicals and unipotent polytopes}, Algebraic Combinatorics, {\bf 1} ( 2018), no. 1, 23-45.

\bibitem{SA}
Andrews S., \textit{Supercharacters of unipotent groups defined by  involutions}, 
J.  Algebra, {\bf 425} (2015),   1-30.


\end{thebibliography}
\end{document}